\documentclass[11pt]{amsart}
\usepackage[english]{babel}
\usepackage[latin1,utf8]{inputenc}
\usepackage{amssymb,geometry,color,xcolor,graphicx}
\usepackage{amsmath,amsthm}
\usepackage{booktabs}
\usepackage{caption}
\usepackage{datetime}
\usepackage{dsfont}
\usepackage{amstext}
\usepackage{fancyhdr}
\usepackage{latexsym}
\usepackage{mathrsfs}
\usepackage{cases}
\usepackage{mathtools}
\usepackage{bm}
\usepackage{subfig}
\usepackage{graphicx}
\usepackage{algorithm}
\usepackage{algpseudocode}
\usepackage{hyperref}
\hypersetup{
    colorlinks=true, 
    linktoc=all,     
    linkcolor=blue,  
}
\numberwithin{equation}{section}

\newcommand{\udef}{\mathrel{\mathop:}=}

\usepackage{dsfont}

\newcommand{\abs}[1]{\left\lvert#1\right\rvert}
\newcommand{\parder}[2]{\frac{\partial#1}{\partial#2}}

\theoremstyle{plain}
\newtheorem{thm}{Theorem}[section]

\newtheorem{defn}[thm]{Definition}
\newtheorem{rem}[thm]{Remark}
\newtheorem{exm}[thm]{Example}

\let\tilde\widetilde
\let\hat\widehat

\renewcommand{\ss}{\scriptstyle}
\def\svdots{\vbox{\baselineskip=1.5pt\lineskiplimit=0pt
	\kern1.5pt \hbox{$\ss .$}\hbox{$\ss .$}\hbox{$\ss .$}}}

\begin{document}

\title[Optimal Control of Moisture Content]{A preliminary model for optimal control of moisture content in unsaturated soils}

\author[Berardi]{Marco Berardi}
\address{Istituto di Ricerca sulle Acque, Consiglio Nazionale delle Ricerche, Via De Blasio, 5, 70132 Bari, Italy}
\email{marco.berardi@ba.irsa.cnr.it}

\author[Difonzo]{Fabio V. Difonzo}
\address{Dipartimento di Matematica, Universit\`a degli Studi di Bari Aldo Moro, Via E. Orabona 4, 70125 Bari, Italy}
\email{fabio.difonzo@uniba.it}

\author[Guglielmi]{Roberto Guglielmi}
\address{Department of Applied Mathematics, University of Waterloo, 200 University Ave W, Waterloo, N2L 3G1 Ontario, Canada}
\email{roberto.guglielmi@uwaterloo.ca}



\null\hfill Version of \today $, \,\,\,$ \xxivtime

\begin{abstract}
In this paper we introduce an optimal control approach to Richards' equation in an irrigation framework, aimed at minimizing water consumption while maximizing root water uptake. We first describe the physics of the nonlinear model under consideration, and then develop the first-order necessary optimality conditions of the associated boundary control problem. We show that our model provides a promising framework to support optimized irrigation strategies, thus facing water scarcity in irrigation. The characterization of the optimal control in terms of a suitable relation with the adjoint state of the optimality conditions is then used to develop numerical simulations on different hydrological settings, that supports the analytical findings of the paper.
\end{abstract}

\maketitle

\section{Introduction}

More and more often, extreme weather events are accompanied by  longer, more intense heat waves and consequent periods of drought, and a forecast global warming is increasing the urgent need of freshwater for human life. In this context, freshwater necessary for agriculture represents almost 70\% of the whole amount of freshwater reserve~\cite{Mao_et_al_COMPAG_2018}.

In this scenario, a wise management of water resources for agricultural purposes is of fundamental importance, even at the irrigation district scale~\cite{Coppola_et_al_water_2019}. Nevertheless, the vast majority of irrigation models just applies heuristic approaches for determining the amount and the timing of irrigation. Albeit an expert knowledge of agricultural and phenological issue is crucial, very seldom such information is coupled with proper mathematical models for controlling irrigation.
In most sophisticated cases,  
the irrigation is managed mainly by studying soil water infiltration into the root zone, starting from Gardner's pioneering works, reviewed in \cite{Gardner1991}; several tools have been proposed in this context, generally providing some type of solver for Richards' equation, the advection-diffusion equation which describes water infiltration in unsaturated porous media accounting also for root water uptake models.\\
In the water resources management or agronomic framework several tools have been proposed in order to benefit from the Richards' equation for irrigation purposes: in  \cite{Difonzo_Masciopinto_Vurro_Berardi} a Python code is presented, able to solve the Richards' equation with any type of root uptake model by a transverse method of lines; also, the popular {\tt Hydrus} software is often used for simulating water flow and root uptake with different crops (e.g apples, as in \cite{NAZARI2021106972},
pecan trees as in \cite{SanjitEtAl2013}).\\

To the best of our knowledge, control methods are seldom applied to irrigation problems, and always with simplified models. For instance, in \cite{Mao_et_al_COMPAG_2018} a zone model predictive control is designed after defining a linear parameter varying model, aimed at maintaining the soil moisture in the root zone within a certain target interval; in \cite{Lopes_optimal_control_irrigation_Hindawi} an optimal control is applied to an irrigation problem, modelled by a simpler (with respect to Richards' equation) hydrologic balance law; a simplified optimization method based on the computation of steady solutions of Richards' equation is proposed in \cite{Berardi_et_al_TiPM_2022}; finally, \cite{Challapa_Molina_IFAC_2019} presents an interesting sliding-mode control approach, but only considering a constant diffusion term in Richards' equation. An elegant approach for applying control techniques in a Richards' equation framework is provided in \cite{Wein_et_al_CNSNS_2019}, yet with very different applications and tools, i.e.  maximizing the amount of absorbed liquid by redistributing the materials, when designing the material properties of a diaper.

In this paper, we propose a model for solving an optimal control problem under the {\em quasi-unsaturated}  assumptions (see Section~\ref{sec:mathmod}), which provide a suitable hydrological setting that prevents to reach water moisture saturation in the soil. We derive the appropriate optimality conditions for the boundary control of a class of nonlinear Richards' equations, and implement these results in the development and computation of numerical solutions by a classical Projected Gradient Descent algorithm. \\
Albeit the focus of this paper does not consist in proposing a novel numerical method, and here a standard {\tt MATLAB} solver is integrated with control tools, some significant advances in the numerical solution of the unsaturated flow model deserve to be reminded; as a matter of fact, the numerics literature on Richards' equation is currently enriching and constantly evolving, since its possibly degenerate and highly nonlinear nature poses several challenges. For instance, the treatment of nonlinearities is a significant issue, and has been faced by different techniques, as Newton methods (\cite{Bergamaschi_Putti,Casulli_Zanolli_SIAM_2010}), L-scheme or its variants \cite{Pop_JCAM_2004,List_Radu_2016,Mitra_Pop_CAMWA_2019}  or Picard iterations \cite{Celia_et_al}.  The discretization in space has been dealt, for instance, by finite elements or mixed finite element methods \cite{Arbogast_Wheeler_1996,Schneid_Knabner_Radu_2004,Kees_Farthing_Dawson_CMAME_2008}, discontinuous Galerkin \cite{Li_Farthing_Miller_ADWR_2007,Clement_et_al_ADWR_2021}, finite volume methods \cite{Eymard_Comput_Geosci_1999,Manzini_Ferraris_ADWR_2004}.
A separate mention is deserved by the problem of infiltration in presence of discontinuities, which can be handled by domain decomposition methods \cite{Seus_Mitra_Pop_Radu_Rohde_2018,Hoang_Pop_Vietnam_2022}, Filippov approach \cite{Berardi_Difonzo_Vurro_Lopez_ADWR_2018,Berardi_Difonzo_Lopez_CAMWA_2020}. For a more detailed discussion, the interested reader is addressed to the following complete reviews on numerical issues in Richards' equation~\cite{Farthing_Ogden_SSSAJ_2017,Zha_et_al_2019_Wires_Water_review}. {The aforementioned methods can be used to face peculiar issues in the numerical integration of Richards' equation, where \texttt{MATLAB pdepe} is known to blow up in a relatively small integration time (see, e.g., \cite{Difonzo_Masciopinto_Vurro_Berardi}). Thus, they could provide further directions in the development of specific algorithms for solving optimal control problems.}

It is worth stressing that there is a vast literature about the problems of existence, uniqueness and regularity of the solution to degenerate parabolic differential equations. In the specific case of Richards' equations, the existence of solutions is tackled in the seminal paper by Alt and Luckhaus \cite{Alt1983}, whose ideas were subsequently developed by other authors. For example,~\cite{Otto1996} describes the semigroup approach to determine the existence of weak solutions, further developed in~\cite{SCHWEIZER2007}. In this paper we adopt the functional framework from~\cite{Marinoschi_book}, that exhaustively describes the minimal assumptions to derive the necessary regularity conditions to develop our analysis, for general classes of hydrological settings. We also refer to \cite{MERZ2010} for a thorough review of such results.

The paper is organized as follows: In Section~\ref{sec:mathmod} we introduce the quasi-unsaturated model of the Richards' equation; in Section~\ref{sec:existence} we describe the framework to ensure the well-posedness of the model of interest, and then derive first order necessary optimality conditions via the Lagrangian method in Section~\ref{sec:optcon}. Finally, we present  numerical simulations in Section~\ref{sec:numericalSimulations}.

\section{The mathematical model}\label{sec:mathmod}
Our work stems from the {\em quasi-unsaturated Richards' model} \cite{Marinoschi_book}, describing a fast diffusion of water in soils. This framework provides a convenient setting to apply optimal control methods and derive optimal irrigation strategies. Indeed, from the mathematical point of view, the quasi-unsaturated diffusive model retains many crucial features of the nonlinear diffusion of water in the soil, while avoiding the special mathematical treatment required to face the limit case of a saturated diffusion.
We consider the model of water diffusion in the space domain $(0,Z)$, where $Z>0$ is the depth of the domain under consideration. We denote by $T\in (0,+\infty)$ the time horizon of the interval $(0, T)$, by $Q = (0,Z)\times (0,T)$ the space-time domain, and by $(\cdot, \cdot)$ and $\|\,\cdot\,\|$ the scalar product and the norm in $L^2(\Omega)$, respectively. In terms of hydraulic parameters, $\theta: Q\to [\theta_r,\theta_S)$ is the \emph{water content} or \emph{moisture}, where $\theta_r$ and $\theta_S$ represent the \emph{residual} and the \emph{saturated water content}, respectively. The function $\beta :[\theta_r,\theta_S)\to [\varrho,+\infty)$ is the \emph{water diffusivity}, satisfying the following condition:
\begin{itemize}
\item[$\mathbf{(H_\beta)}$] $\beta$ is locally Lipschitz  continuous and monotonically increasing, $\beta(\theta) \ge \varrho  > 0$ for all $\theta \in [\theta_r, \theta_S)$, and $\displaystyle\lim_{\theta\nearrow \theta_S}
\beta(\theta) = +\infty$.
\end{itemize}

\noindent
The function $\beta^*$ is the primitive of the water diffusivity $\beta$ that vanishes \mbox{at $0$.}
Thus, assumption $(\mathbf{H}_\beta)$ implies that $\beta^*$ is differentiable and monotonically increasing on $[\theta_r, \theta_S)$, and satisfies
$$
(\beta^*(\theta_1) - \beta^*(\theta_2))(\theta_1 - \theta_2) \ge \varrho (\theta_1 - \theta_2)^2\qquad\forall \;\theta_1, \theta_2 \in [\theta_r, \theta_S)\, .
$$

\noindent
Moreover, the \emph{hydraulic conductivity} $K:[\theta_r,\theta_S]\to \mathbb{R}$ is a non-negative, Lipschitz continuous on $[\theta_r, \theta_S]$, and monotonically increasing function.
Other hydraulic functions of interest are the \emph{{liquid} pressure head $h:Q\to (-\infty,0)$}, {that is negative for unsaturated porous media}, and the \emph{specific water capacity} $C(h) = \frac{\mathrm{d}\theta}{\mathrm{d}h}$, {that practically represents a storage term}. The relation between the functions $\beta$, $K$ and $C$ is then expressed by \mbox{$\beta(\theta(h)) := \frac{K(h)}{C(h)}$}. \\
With these notations, and assuming the vertical axis with downward positive orientation, the implicit form of the quasi-unsaturated model of the nonhysteretic infiltration of an incompressible fluid into an isotropic, homogeneous, unsaturated porous medium with a constant porosity and truncated diffusivity with non-homogeneous--Dirichlet Boundary Conditions (BCs) is given by the system
\begin{equation}\label{eq:1DRichImplicitTrunc}
\begin{cases}
\frac{\partial\theta}{\partial t} - \frac{\partial^2\beta^*(\theta)}{\partial z^2} + \frac{\partial K(\theta)}{\partial z} = f(\theta) &
 \text{ in }\ Q \, ,\\
\theta(z,0) = \theta_0(z) &
 \text{ in }\ (0,Z)\, , \\
\theta(0,t) = v(t) & \text{ for }\ t\in (0,T)\, , \\
\theta(Z,t) = g(t) & \text{ for }\ t\in (0,T)\, .
\end{cases}
\end{equation}

The source term $f(\theta)$ represents a sink function, that in our model describes the root water uptake.
The non-homogeneous Dirichlet BCs are given by two functions $v,g:(0,T)\to \mathbb{R}$, where the BCs at $z = 0$ is the control input $v$, that describes the irrigation strategy over the time horizon $(0,T)$, while the BCs at $z = Z$ is a given function $g$ modeling the interaction with the environment below the root zone
\begin{equation}\label{eq:constraints}
\theta_r <
v(t),g(t) < \theta_S \quad\text{ for a.e. }\ t\in (0,T)\;.
\end{equation}

\begin{rem}
Let us notice that the first equation in \eqref{eq:1DRichImplicitTrunc} is the Richards' equation. In fact, we can easily compute that
\begin{align*}
\parder{}{z}\beta^* &= \parder{\beta^*}{\theta}\parder{\theta}{h}\parder{h}{z} = \frac{K(h)}{C(h)}C(h)\parder{h}{z} = K(h)\parder{h}{z},
\end{align*}
that provides the diffusion term in the classical mixed form of the Richards' equa\-tion~\cite[Eq. (3)]{Celia_et_al}.
\end{rem}

\section{Existence of solutions}\label{sec:existence}

Following~\cite[Section~4.5]{Marinoschi_book}, we first reduce system~\eqref{eq:1DRichImplicitTrunc} to a problem with homogeneous Dirichlet BCs. For this purpose, we assume the following hypothesis
\begin{align*}
&\mathbf{(H_\omega)}\quad \exists\; \omega:Q\to\mathbb{R}\ \text{ such that } \\
&\begin{cases}
\omega\in L^2(0,T; H^{1}(0,Z))\cap L^\infty(Q)\, ,\quad
\omega_t\in L^2(Q) \,, & \\
\|\omega\|_{L^\infty(Q)} < \theta_S\, ,\quad \mathrm{essinf}_{(z,t)\in Q}\; \omega(z,t) >\theta_r\; ,& \\
\omega(0,t) = v(t)\,,\ \omega(Z,t) = g(t)\,,\ \text{ for a.e. }\ t \in (0,T)\, , &
\end{cases}
\end{align*}
where $\omega_t$ denotes the time derivative in the sense of distributions from $(0, T)$ to $L^2(0,Z)$. Thus, the function $\omega$ is defined on the cylinder $Q$ and it attains the boundary values of $v$ and $g$ in $z=0$ and $z = Z$, respectively. Then, introducing the function $\phi = \theta - \omega$, 
system~\eqref{eq:1DRichImplicitTrunc} is equivalent to
\begin{equation}\label{eq:equiv1DRichTrunc}
\left\{
{
\begin{array}{ll}
\frac{\partial\phi}{\partial t} - \frac{\partial^2 F^\omega(\phi)}{\partial z^2} + \frac{\partial K(\phi + \omega)}{\partial z} = f^B(\phi) - \omega_t &
 \text{ in }\ Q \, ,\\[.8ex]
\phi(z,0) = \phi_0(z) &
 \text{ in }\ (0,Z)\, , \\[0.4ex]
\phi(0,t) = 0 & \ t\in (0,T)\, , \\[0.4ex]
\phi(Z,t) = 0 & \ t\in (0,T)\, ,
\end{array}
}
\right.
\end{equation}
where $\phi_0 := \theta_0 - \omega_0$ and, for all $\phi\in V$,
\[
F^\omega(\phi) := \beta^*(\phi + \omega) - \beta^*(\omega)\; ,\qquad
f^B(\phi) := f(\phi + \omega) + \frac{\partial^2\beta^*(\omega)}{\partial z^2}\; .
\]

We now introduce a suitable functional framework for~\eqref{eq:equiv1DRichTrunc}. Let us consider the Gelfand triple $H = L^2(0,Z)$, $V = H^1_0(0,Z)$, and its dual $V' = H^{-1}(0,Z)$, with their usual norms.
Then,~\eqref{eq:equiv1DRichTrunc} is equivalent to the abstract equation
\begin{equation}\label{eq:equiv1Trunc}
\frac{\mathrm{d}\phi}{\mathrm{d} t} + B(t)\phi = f^B(\phi) - \omega_t\quad \text{ a.e. }
\ t\in (0,T) \, ,\qquad
\phi(0) = \phi_0\, ,
\end{equation}
where the operator $B(t):V\to V'$ is defined by
\[
\langle B(t)\phi,\psi\rangle_{V',V} = \int_0^Z \left(\frac{\partial F^\omega(\phi(t))}{\partial z} - K(\phi(t) + \omega(t))\right)\frac{\partial\psi}{\partial z}\mathrm{d}z\qquad \forall \phi, \psi\in V\;.
\]
Notice that, thanks to $\mathbf{(H_\omega)}$, we have that $\frac{\partial^2\beta^*(\omega)}{\partial z^2}\in L^2(0,T;V')$. Indeed, for all $\varphi\in V = H^1_0(0,Z)$,
\[
\abs{\int_0^Z
-\frac{\partial^2 \beta^*(\omega)}{\partial z^2} \varphi \,\mathrm{d}z
} =
\abs{
\int_0^Z
\frac{\partial \beta^*(\omega)}{\partial z} \frac{\partial\varphi}{\partial z} \mathrm{d}z
}
\le
\left\|\beta(\omega)\frac{\partial\omega}{\partial z}\right \|\, \left\|\frac{\partial\varphi}{\partial z}\right \|\le M_\omega \|\varphi\|_V\, .
\]
Hence $\|\frac{\partial^2 \beta^*(\omega)}{\partial z^2}\|_{V'}\le M_\omega$, and thus $f^B - \omega_t\in L^2(0,T;V')$. In particular, this implies that the right-hand side of~\eqref{eq:equiv1Trunc} is also in $ L^2(0,T;V')$.
\begin{defn}
Let $\theta_0\in L^2(0,Z)$ and $f:[\theta_r,\theta_S]\to\mathbb{R}$ be Lipschitz continuous. We say that the function $\phi\in C([0,T];L^2(0,Z))$ is a solution to~\eqref{eq:equiv1Trunc} if $\frac{\mathrm{d}\phi}{\mathrm{d}t}\in L^2(0,T;V')$, $F^\omega\in L^2(0,T;V)$, $\phi(0) = \phi_0$ and, for a.e. $t\in (0,T)$ and for all $\psi\in V$,
\begin{equation}\label{eq:def1}
\left\langle \frac{\mathrm{d}\phi}{\mathrm{d}t}(t),\psi \right\rangle_{V',V} + \langle B(t)\phi,\psi\rangle_{V',V} = \left\langle f^B(\phi) - \omega_t,\psi\right\rangle_{V',V} \; .
\end{equation}
\end{defn}
Existence of solutions to~\eqref{eq:equiv1Trunc} with a source term independent of the water content $\theta$ is proved in~\cite[Proposition~5.3]{Marinoschi_book}. For our purpose, we need to extend such well-posedness result to the case of a nonlinear source term $f(\theta)$, as in system~\eqref{eq:1DRichImplicitTrunc}. This can indeed be achieved by following a standard Galerkin approximation approach (see, e.g.,~\cite[Lemma~5.3]{troltzsch2010optimal}).
\begin{thm}\label{thm:wellposed}
Assume $\mathbf{(H_\omega)}$, $\theta_0 \in L^2(0,Z)$, and that $f:[\theta_r,\theta_S]\to\mathbb{R}$ is Lipschitz continuous. Then the problem~\eqref{eq:equiv1Trunc} admits a unique solution $\phi \in C([0, T];L^2(0,Z)) \cap L^2(0, T;V)$ with $\frac{\mathrm{d}\phi}{\mathrm{d}t}\in L^2(0, T; V')$ and $
F^\omega(\phi)\in L^2(0, T; V)$.
Therefore, system~\eqref{eq:1DRichImplicitTrunc} admits a unique solution $\theta \in L^2(0,T;H^1(0,Z))\cap C([0,T];H)$ with $\frac{\mathrm{d}\theta}{\mathrm{d}t} \in L^2(0,T;V')$ and $\beta^*(\theta)\in L^2(0,T;H^1(0,Z))$.
\end{thm}
Moreover, for appropriate initial conditions, we can prove that the solution stays away from the saturation value $\theta_S$ uniformly in time. To this aim, we introduce the function $j:\mathbb{R}\to (-\infty,\infty]$ defined by
\[
j(r)\udef
\begin{cases}
\displaystyle \int_0^r \beta^*(\xi) d\xi\, , & r < \theta_S\, ,\\[0.5ex]
\displaystyle \infty\, , & r \ge \theta_S\, ,
\end{cases}
\]
and the space
\[
M_j \udef \{\theta\in L^2(0,Z) : j(\theta)\in L^1(0,Z)\}\, .
\]
\begin{thm}[Theorem 5.7,~\cite{Marinoschi_book}]\label{thm:bndaway}
Assume $\mathbf{(H_\omega)}$, $\theta_0 \in M_j$, and that $f:[\theta_r,\theta_S]\to\mathbb{R}$ is Lipschitz continuous. Moreover, assume that $f$ is non-negative, that is, there exists $f_m\in [0,\infty)$ such that $f_m\le f$, and
\begin{eqnarray*}
& \mathrm{essinf}_{x\in (0,Z)} \theta_0(x)\ge 0\, ,\\[0.5ex]
& \theta_m(t)\le g(t),v(t) < \theta_S\, ,\quad \text{for all }\ t\in [0,T]\, ,
\end{eqnarray*}
where
\[
\theta_m(t) = \mathrm{essinf}_{x\in (0,Z)} \theta_0(x) + f_m t\, .
\]
Then the solution $\theta$ to problem~\eqref{eq:1DRichImplicitTrunc} satisfies
\[
\theta_m(t)\le \theta(x,t) < \theta_S\, ,\quad \text{for all }\ (x,t)\in \{0,Z\}\times [0,T]\, .
\]
\end{thm}

\section{The optimal control problem}\label{sec:optcon}
\noindent
In this section we formally derive the first order necessary optimality conditions for the cost functional
\begin{equation}\label{eq:costfun}
J(\theta,u) = \frac 12 \int_Q \abs{f(\theta(z,t)) - 1}^2 \,\mathrm{d}z\mathrm{d}t + \frac{\lambda}{2} \int_0^T \abs{u(t)}^2\,\mathrm{d}t\; ,
\end{equation}
where $u = v - \theta_r$, $v$ is the control that appears in~\eqref{eq:1DRichImplicitTrunc}, $\lambda > 0$ is the coefficient of the control cost, $f:[\theta_r,\theta_S]\to\mathbb{R}$ describes the normalized root water uptake model as in \eqref{eq:1DRichImplicitTrunc}, and $\theta$ is the solution to~\eqref{eq:1DRichImplicitTrunc} with $f$ as the sink term.
Roughly speaking, the performance index~\eqref{eq:costfun} optimizes the root water uptake (see, for example, expression~\eqref{eq:Feddes} in Section \ref{sec:numericalSimulations}, where $f$ is maximized when $f\equiv1$) while minimizing the irrigation cost $u$. In this setting, it is natural to consider the following space of admissible control
\begin{equation}\label{eq:inputconstraints}
U_{ad} := \{u\in L^\infty(0,T) : 0\le u(t) < \theta_S - \theta_r\ \text{ for a.e. }\ t\in (0,T)\}\;.
\end{equation}
{Fixing $g \in L^2(0, T)$, $\theta_0 \in L^2(0,Z)$, we introduce the control-to-state operator $\Lambda: U_{ad} \to C([0, T ];L^2(0,Z))$ such that $u\in U_{ad} \mapsto \theta \in C([0, T ]; L^2(0,Z))$ solution of~\eqref{eq:1DRichImplicitTrunc}.}
Theorem~\ref{thm:wellposed} ensures that the mapping $\Lambda$ is well-posed. We can thus reformulate the minimization of a functional $\tilde{J}(\theta,u)$ constrained to the control system~\eqref{eq:1DRichImplicitTrunc} in terms of the so-called \emph{reduced cost functional} $J : U_{ad}\to \mathbb{R}$ defined by $J (u) := \tilde{J}(\Lambda(u), u)$.
We first introduce the Lagrangian functional
\begin{align*}
\mathcal{L}(\theta,u,{\bf p}) &= J(\theta,u)
- \int_Q \left[\frac{\partial\theta}{\partial t} - \frac{\partial}{\partial z}\left(\beta\frac{\partial\theta}{\partial z}\right) + \frac{\partial K(\theta)}{\partial z} - f\right] p\,\mathrm{d}z\mathrm{d}t \\
&\quad - \int_0^T\left(\theta(0,t) - u(t)\right)p_1\,\mathrm{d}t - \int_0^T\left(\theta(Z,t) - g(t)\right)p_2\,\mathrm{d}t,
\end{align*}
where ${\bf p} = (p,p_1,p_2)$ are adjoint variables that will be useful to find a representation of the optimal control. After integration by parts, we can rewrite the Lagrangian functional as
\begin{align*}
\mathcal{L}(\theta,u,{\bf p}) &= J(\theta,u)
- \int_Q \frac{\partial\theta}{\partial t} p + \beta\frac{\partial\theta}{\partial z}\frac{\partial p}{\partial z} + \left( \frac{\partial K(\theta)}{\partial z} - f\right)p\,\mathrm{d}z\mathrm{d}t \\
&\quad + \int_0^T\left[\left(\beta \frac{\partial\theta}{\partial z} p\right)_{\big\vert z = Z} - \left(\beta \frac{\partial\theta}{\partial z} p\right)_{\big\vert z = 0}\right]\,\mathrm{d}t +\\
&\quad - \int_0^T\left(\theta(0,t) - u(t)\right)p_1\,\mathrm{d}t - \int_0^T\left(\theta(Z,t) - g(t)\right)p_2\,\mathrm{d}t\; .
\end{align*}
Hereafter, we shall assume that the source term $f\in H^1(\theta_r,\theta_S)$ to justify the following computations.
In order to derive the first order optimality conditions of problem~\eqref{eq:costfun}-\eqref{eq:1DRichImplicitTrunc} with input constraints~\eqref{eq:inputconstraints}, we enforce the condition $D_\theta\mathcal{L}(\theta^*,u^*,p^*)\theta = 0$ for all $\theta$, that determines the equation satisfied by the adjoint variable $p$; and the condition $D_u\mathcal{L}(\theta^*,u^*,p^*)\cdot (u - u^*)\ge 0$ for all $u\in U_{ad}$, that returns the optimality condition satisfied by any optimal control $u^*$. After direct computations, we get that
\begin{multline*}
D_\theta \mathcal{L}(\theta^*,u^*,{\bf p})\theta =  - \int_0^Z\left[\theta(z,T) p(z,T) - \theta(z,0) p(z,0)\right]\,\mathrm{d}z\\
\quad + \int_Q \theta \left[\frac{\partial p}{\partial t} + (1 - f(\theta^*)) \frac{\mathrm{d}f}{\mathrm{d}\theta}(\theta^*) + \beta\frac{\partial^2 p}{\partial z^2} + \frac{\mathrm{d}K}{\mathrm{d}\theta}(\theta^*) \frac{\partial p}{\partial z}\right] \,\mathrm{d}z\mathrm{d}t \\
\quad + \int_0^T\left[\frac{\partial\theta}{\partial z}(Z,t) \left(\beta(\theta^*)p\right)_{\big\vert z = Z} - \frac{\partial\theta}{\partial z}(0,t) \left(\beta(\theta^*)p\right)_{\big\vert z = 0}\right]\,\mathrm{d}t\\
\quad + \int_0^T \theta(Z,t)\left[- \beta(\theta^*) \frac{\partial p}{\partial z} - \frac{\mathrm{d}K}{\mathrm{d}\theta}(\theta^*) p + \frac{\mathrm{d}\beta}{\mathrm{d}\theta}(\theta^*) \frac{\partial\theta^*}{\partial z} p - p_2\right]\,\mathrm{d}t\\
\quad - \int_0^T \theta(0,t)\left[- \beta(\theta^*) \frac{\partial p}{\partial z} - \frac{\mathrm{d}K}{\mathrm{d}\theta}(\theta^*) p + \frac{\mathrm{d}\beta}{\mathrm{d}\theta}(\theta^*)  \frac{\partial\theta^*}{\partial z} p + p_1\right]\,\mathrm{d}t \; .
\end{multline*}
Thus, we deduce that the adjoint variable ${\bf p}$ satisfies
\[
\begin{cases}
\frac{\partial p}{\partial t} + \beta(\theta^*) \frac{\partial^2 p}{\partial z^2} + \frac{\mathrm{d}K}{\mathrm{d}\theta}(\theta^*) \frac{\partial p}{\partial z} = F(\theta^*) &
 \text{ in } Q,\\
p(z,T) = 0 &
 \text{ in }\ (0,Z)\, , \\
p(0,t) = p(Z,t) = 0 & \text{ for }t\in (0,T)\, , \\
p_1 = \left(\beta(\theta^*) \frac{\partial p}{\partial z} \right)_{\big\vert z = 0} & \\
p_2 = -\left(\beta(\theta^*) \frac{\partial p}{\partial z}
\right)_{\big\vert z = Z} &
\end{cases}
\]

where
\[
F(\theta)\udef[f( \theta) - 1] \frac{\mathrm{d}f}{\mathrm{d}\theta}(\theta).
\]
On the other hand, since
\begin{align*}
D_u \mathcal{L}(\theta^*,u^*,{\bf p})u &=
\int_0^T (\lambda u^* + p_1) u\,\mathrm{d}t\; ,
\end{align*}
the condition $D_u\mathcal{L}(\theta^*,u^*,p^*)\cdot (u - u^*)\ge 0$ for all $u\in U_{ad}$ implies the optimality condition
\[
\left\langle\lambda u^*(t) + \left(\beta(\theta^*) \frac{\partial p}{\partial z}
\right)_{\big\vert z = 0},{u} - u^*\right\rangle_{L^2(0,T)} \ge 0
\]
for all ${u}\in U_{ad}$.
We thus obtain that any optimal solution $(\theta^*,u^*,p^*)$ of problem~\eqref{eq:1DRichImplicitTrunc}-\eqref{eq:costfun}-\eqref{eq:inputconstraints} must satisfy the optimality system
\begin{subequations}\label{eq:optSys}
\begin{align}
&\begin{cases}
\frac{\partial\theta}{\partial t} - \frac{\partial^2\beta^*(\theta)}{\partial z^2} + \frac{\partial K(\theta)}{\partial z} = f &
 \text{ in }Q,\\
\theta(z,0) = \theta_0(z) &
 \text{ in }(0,Z), \\
\theta(0,t) = v(t) & \text{ for }t\in (0,T), \\
\theta(Z,t) = g(t) & \text{ for }t\in (0,T),\\
\end{cases}
\\
&
\begin{cases}
\frac{\partial p}{\partial t} + \beta(\theta^*) \frac{\partial^2 p}{\partial z^2} + \frac{\mathrm{d}K}{\mathrm{d}\theta}(\theta^*) \frac{\partial p}{\partial z} = F(\theta^*) &
 \text{ in } Q,\\
p(z,T) = 0 &
 \text{ in }\ (0,Z)\, , \\
p(0,t) = p(Z,t) = 0 & \text{ for }t\in (0,T)\, ,
\end{cases} \\[2ex]
& \left\langle\lambda u^*(t) + \left(\beta(\theta^*) \frac{\partial p}{\partial z}
\right)_{\big\vert z = 0},{u} - u^*\right\rangle_{L^2(0,T)} \ge 0
\end{align}
\end{subequations}
 for all ${u}\in U_{ad}$, where we recall that $v = \theta_r + u$. In the next section, we exploit this optimality system to build suitable algorithms to numerically solve the optimal control problem~\eqref{eq:costfun}-\eqref{eq:1DRichImplicitTrunc}.

\section{Algorithm and numerical simulations}\label{sec:numericalSimulations}

\begin{figure}
\centering
\includegraphics[width=0.8\linewidth]{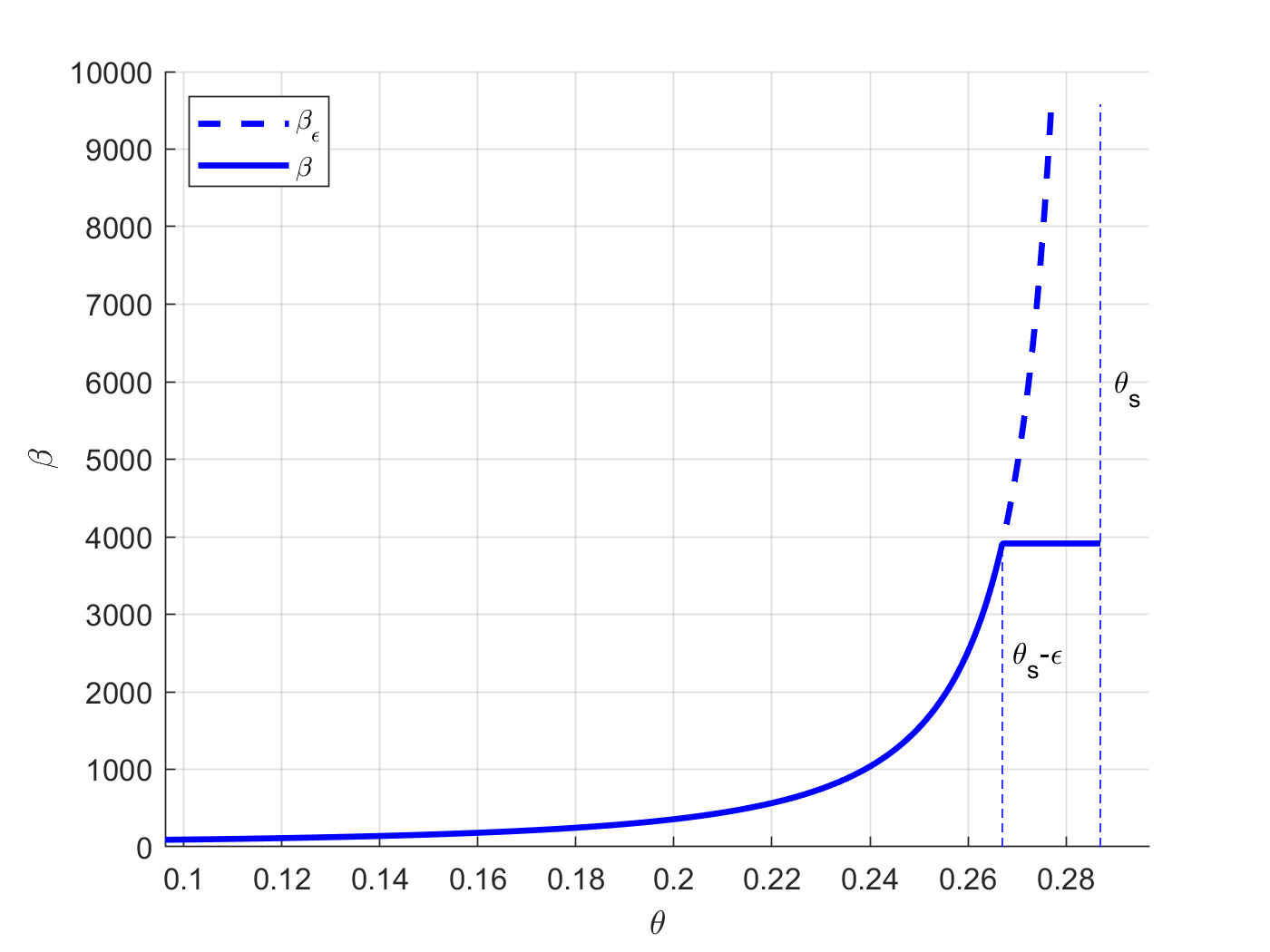}
\caption{Graph of the truncated water diffusivity $\beta_\varepsilon$ from \eqref{eq:betaEps}.}
\label{fig:betaEps}
\end{figure}

Our optimization procedure will follow the Projected Gradient Descent (PGD) described in Algorithm \ref{alg:PGD} (see \cite{troltzsch2010optimal} for a thorough introduction to such optimization algorithms).
However, when solving \eqref{eq:optSys} with PGD, it could happen that Theorem \ref{thm:bndaway} is not satisfied at each iteration, thus incurring numerical difficulties due to the singularity of water diffusivity $\beta$ at $\theta = \theta_S$. Therefore, we shall approximate it by truncation: given a small $\varepsilon > 0$, we define
\begin{equation}\label{eq:betaEps}
\beta_\varepsilon(r)\udef
\begin{cases}
\beta(r), & r\le \theta_S - \varepsilon, \\
\beta(\theta_S -\varepsilon), & r > \theta_S -\varepsilon,
\end{cases}
\end{equation}
as shown in the Figure \ref{fig:betaEps}.
Regularization \eqref{eq:betaEps} is a standard technique when dealing with Richards' equation to handle singularities in the diffusion term, and it is used in finite difference schemes \cite{Alt1983,Pop2011RegularizationSF} or FEM \cite{Nochetto1997} for both the mathematical and numerical analysis of degenerate, and possibly doubly-degenerate, parabolic equations. \\
In the following simulations, we are then actually computing the numerical solutions to \eqref{eq:optSys} after replacing $\beta$ with $\beta_\varepsilon$ as defined in \eqref{eq:betaEps}. \\
Moreover, we have set a maximum number of iterations equal to 100 before exiting PGD iterations, a tolerance of $10^{-5}$, and a regularization parameter $\varepsilon=10^{-3}$ for water diffusivity in \eqref{eq:betaEps}. We stress that the order of magnitude of $\varepsilon$ has been chosen so to be consistent with that of the different $\theta_S$ values selected in all the simulations that follow.

\begin{algorithm}
\caption{Projected Gradient Descent Algorithm applied to Richards' equation for determining optimal irrigation. Intermediate steps for solving direct (line 3) and adjoint (line 4) problems and for computing the optimal descent direction step-size (line 6) are performed using \texttt{MATLAB} \texttt{pdepe} and \texttt{fmincon} functions, respectively.}
\label{alg:PGD}
\begin{algorithmic}[1]
\State $n\gets1, u_n\gets u_0$\Comment{Initialization}
\While{$n\leq\textrm{maxit}$}
\State $\theta_n\gets$ \Call{directProblem}{$u_n$}\Comment{Solving the I-BV problem
}
\State $p_n\gets$ \Call{adjointProblem}{$\theta_n$}\Comment{Solving the adjoint problem
}
\State $r_{n}\gets-\left(\lambda u_{n}+\beta(p_{n}(z=0))\partial_z p_n(z=0)\right)$\Comment{Descent direction}
\State $s_{n}\gets\arg\min_{s} J(\mathrm{pr}_{\mathbf{U}}(u_n+s r_n))$\Comment{Optimal direction Step-size
}
\Statex
\If{$\abs{J(\mathrm{pr}_{\mathbf{U}}(u_n+s_n r_n))-J(u)}<\mathrm{tol}$}\Comment{Exit condition}
\State \textbf{break}
\Else
\State $u\gets\mathrm{pr}_{\mathbf{U}}(u_n+s_n r_n)$\Comment{Updating control}
\State $n\gets n+1$
\EndIf
\EndWhile
\end{algorithmic}
\end{algorithm}

Moreover, we select a root water uptake model of Feddes type (used, for instance, in \cite{Difonzo_Masciopinto_Vurro_Berardi,SanjitEtAl2013}) as source term in \eqref{eq:1DRichImplicitTrunc}. Its expression is given by
\begin{equation}\label{eq:Feddes}
f(h)=\varphi\hat f(h),\quad
\hat f(h)\udef
\begin{cases}
0, & \textrm{ if }h_{1}\leq h\leq0\textrm{ or } h\leq h_{4}, \\
\frac{ h-h_{1}}{h_{2}-h_{1}}, & \textrm{if }h_{2}< h<h_{1}, \\
1, & \textrm{if\quad}h_{3}\leq  h \leq h_{2}, \\
\frac{ h-h_{4}}{h_{3}-h_{4}}, & \textrm{if }h_{4}< h<h_{3},
\end{cases}
\end{equation}
with the following values, in cm: $h_4\approx-820,\ h_3\approx-400,\  h_2\approx-350,\ h_1=0$. Also, we set
$\varphi=0.1/Z$, where $Z$ is the soil depth, and $\lambda=0.1$ in \eqref{eq:costfun}.

\begin{rem}
Let us notice that the maximum value for $\hat f(h)$ in \eqref{eq:Feddes} is set to $1$ for normalization purposes. In fact, when it comes to practical problems, one uses $f(h)$ properly rescaled according to experimental evidences through the factor $\varphi$, which is the ratio of the potential transpiration rate and the rooting depth, as explained in \cite{Utset2000}. \\
Moreover, we stress that in general one does not necessarily require the source term $f(h)$ to be zero for the values of $h$ corresponding to the boundary of $[\theta_r,\theta_S]$. However, from a physical point of view, it makes sense for a source term to vanish when the soil is either dry or saturated. This is exactly the case of Feddes-type source terms as the one we consider in \eqref{eq:Feddes} for our numerical simulations.
\end{rem}

In Example \ref{ex:haverkamp1} and Example \ref{ex:haverkamp2} below, we simulate a soil described by Haverkamp model \cite{Celia_et_al,Berardi_Difonzo_Notarnicola_Vurro_APNUM_2019}, whose constitutive relations are given by
\begin{equation}\label{eqn:Haverkamp}
\theta( h) = \frac{\alpha \left(\theta_S - \theta_r \right)}{\alpha + \abs{h}^{\beta_2}} + \theta_r\; ,\qquad
K( h) = K_S  \frac{A}{A+\abs{h}^{\beta_1}}\; ,
\end{equation}
representing water retention curve and hydraulic conductivity, respectively. We first verify that Haverkamp model falls within the quasi-unsaturated model for a suitable choice of the parameters involved in the model setting.
In fact, in Haverkamp model we have
\begin{equation}\label{eq:HaverModel}
\beta(\theta(h))=\frac{K_S A(\alpha+(-h)^{\beta_2})^2}{(A+(-h)^{\beta_1})\alpha(\theta_S-\theta_r)\beta_2(-h)^{\beta_2-1}}
\end{equation}
with $\beta_1,\beta_2>0$, that is Lipschitz, monotonically increasing and bounded from below. In order to satisfy assumption ($\mathbf{H_\beta}$), we shall ensure that
\begin{equation*}
\lim_{h\nearrow0^-}\beta(\theta(h))=+\infty.
\end{equation*}
From~\eqref{eq:HaverModel}, a straightforward computation provides that this condition is satisfied if and only if $\beta_2>1$.

\begin{exm}\label{ex:haverkamp1}
This is the case of the sandy soil considered in \cite{Celia_et_al}, with parameters
\begin{equation}\label{eq:parametersEx}
\begin{split}
K_S=34\,{\rm cm/h},\,\, A=1.175\times 10^6,\,\,\beta_1=4.74,\\
\theta_S=0.287,\,\,\theta_r=0.075,
\alpha=1.611\times10^6,\,\,\beta_2=3.96,
\end{split}
\end{equation}
where the root water uptake model is as in \eqref{eq:Feddes}.

We have performed a simulation of $T=3$ hours with a maximum depth of $Z=70$ cm.

As in~\eqref{eq:1DRichImplicitTrunc}, boundary condition at the top varies in time according to the irrigation strategy:
\begin{equation}\label{eq:topBC_ex1}
\theta(0,t)=\theta_{\rm top}(t)=u(t)+\theta_r,
\end{equation}
while bottom condition has been chosen so to be constant over time:
\begin{equation}\label{eq:bottomBC_ex1}
\theta(Z,t)=\theta_{\rm bottom}(t)=0.9\theta_r+0.1\theta_S,\,\,t\in[0,T].
\end{equation}
Finally, initial condition is linearly varying over time as
\begin{equation}\label{eq:IC_ex1}
\theta(z,0)=\theta_{\rm top}(0)+z\frac{\theta_{\rm bottom}(0)-\theta_{\rm top}(0)}{Z},\quad z\in[0,Z].
\end{equation}

\begin{figure}
 \centering
 \subfloat[][Optimal water content in the spatio-temporal domain.]
 {\includegraphics[width=.48\textwidth]{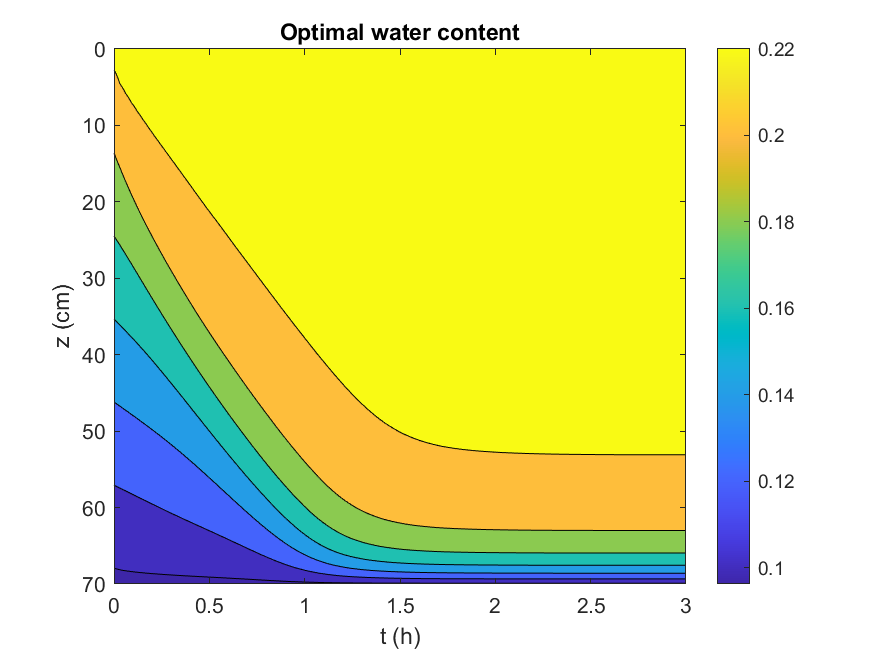}\label{fig:optTheta1}} \quad
 \subfloat[][Optimal water content profiles over time.]
 {\includegraphics[width=.48\textwidth]{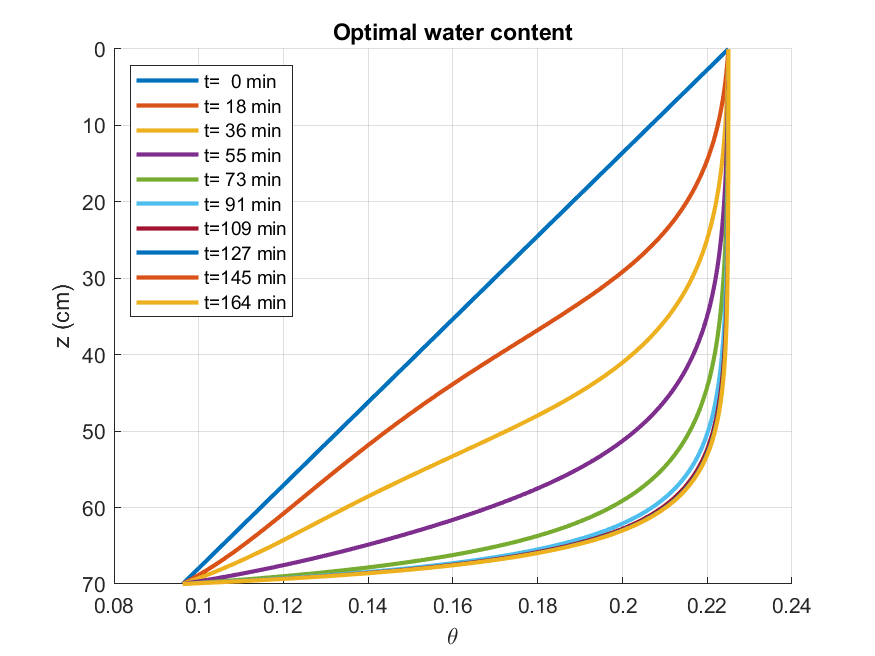}\label{fig:optMeanTheta1}} \\
 \subfloat[][Average optimal water content and optimal control.]
 {\includegraphics[width=.48\textwidth]{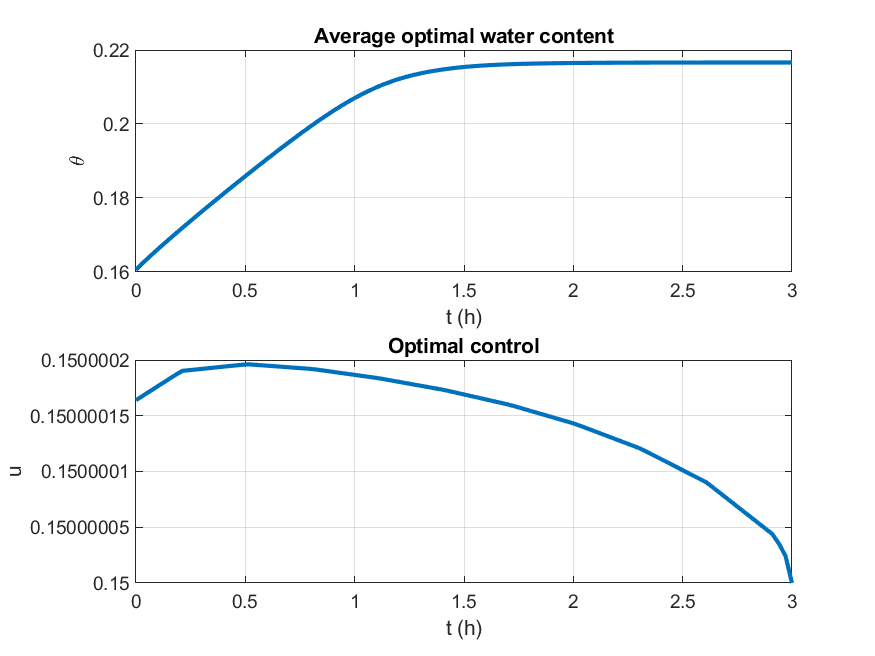}\label{fig:optControl1}}
 \caption{Numerical simulations relative to Example \ref{ex:haverkamp1}, where initial condition is given by \eqref{eq:IC_ex1} and boundary conditions are as in \eqref{eq:topBC_ex1} and \eqref{eq:bottomBC_ex1}, respectively.
 }
  \label{fig:optResults1}
 \end{figure}

Our simulations have been produced using \texttt{MATLAB} \texttt{active-set} algorithm; we report that same results are obtained using \texttt{MATLAB sqp}. We have observed that convergence is reached after 3 iterates within the given tolerance, and the numerical solution is locally optimal.
Results are in Figure~\ref{fig:optResults1}. As can be seen, the optimal control framework succeeds in determining an optimal control that optimizes the performance index~\eqref{eq:costfun}, with a reduced water consumption and average water content over time.
\end{exm}

\begin{exm}\label{ex:haverkamp2}
In this second simulation, using the same soil as in Example \ref{ex:haverkamp1}, we consider a time-varying bottom condition
\begin{equation}\label{eq:bottomBC_ex2}
\theta(Z,t)=\theta_{\textrm{bottom}}(t)=\left(1-\frac{t}{T}\right)\theta_{b1}+\frac{t}{T}\theta_{b2},\,\,t\in[0,T],
\end{equation}
where
\[
\theta_{b1}\udef0.9\theta_r+0.1\theta_S,\,\theta_{b2}\udef0.7\theta_r+0.3\theta_S,
\]
while top condition is given by
\begin{equation}\label{eq:topBC_ex2}
\theta(0,t)=\theta_\textrm{top}(t)=u(t)+\theta_r.
\end{equation}
Moreover, initial condition is
\begin{equation}\label{eq:IC_ex2}
\theta(z,0)=\theta_{\rm top}(0)+z\frac{\theta_{\rm top}(0)-\theta_{\rm bottom}(0)}{Z},\quad z\in[0,Z].
\end{equation}
It turns out that \texttt{MATLAB sqp} converges, in $3$ iterates, to a local optimal solution, further providing the best results if compared to \texttt{active-set}. Results are displayed in Figure \ref{fig:optResults2}.
\begin{figure}
 \centering
 \subfloat[][Optimal water content in the spatio-temporal domain.]
 {\includegraphics[width=.48\textwidth]{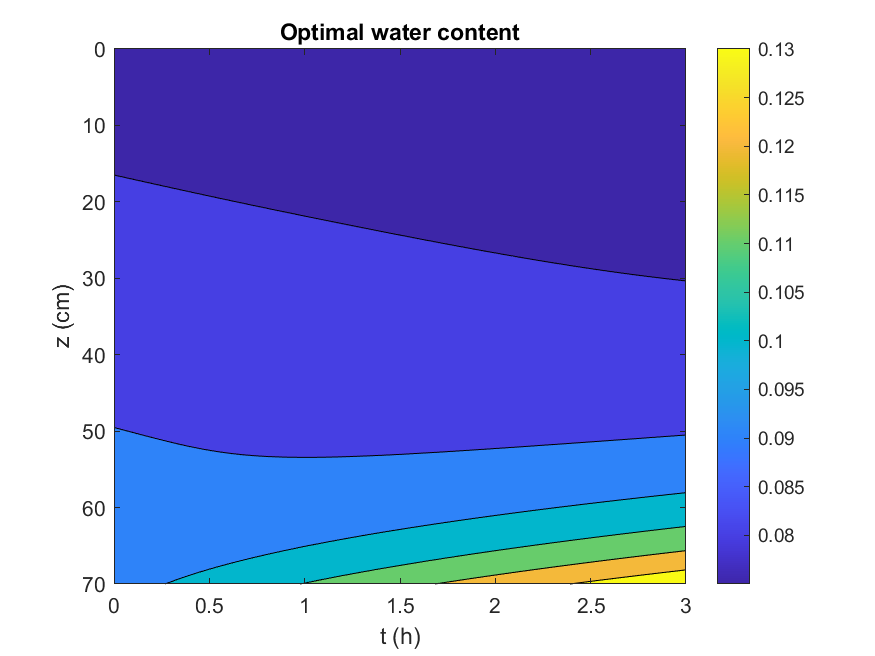}\label{fig:optTheta2}} \quad
 \subfloat[][Average optimal water content over time.]
 {\includegraphics[width=.48\textwidth]{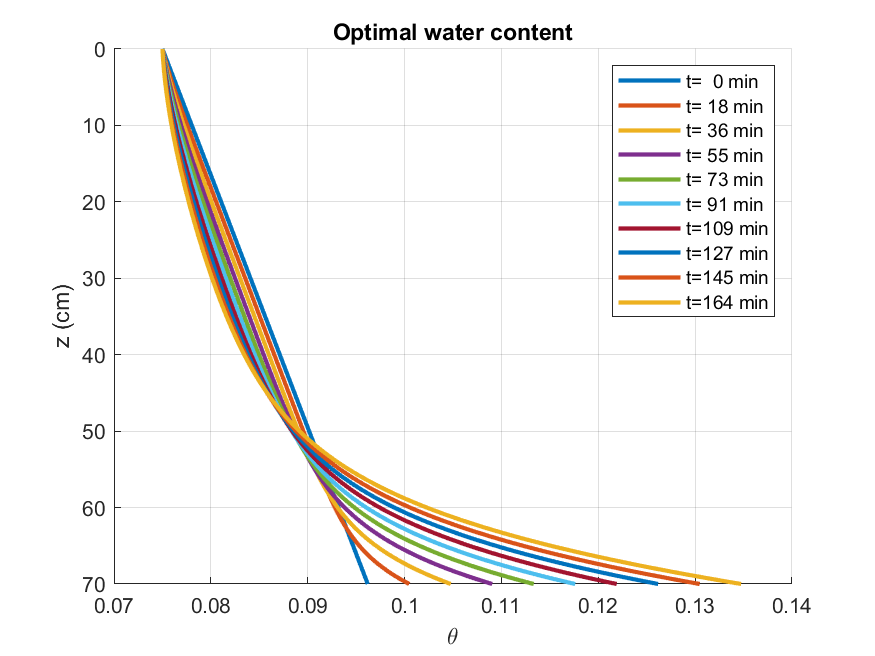}\label{fig:optMeanTheta2}} \\
 \subfloat[][Optimal water content profiles and optimal control.]
 {\includegraphics[width=.48\textwidth]{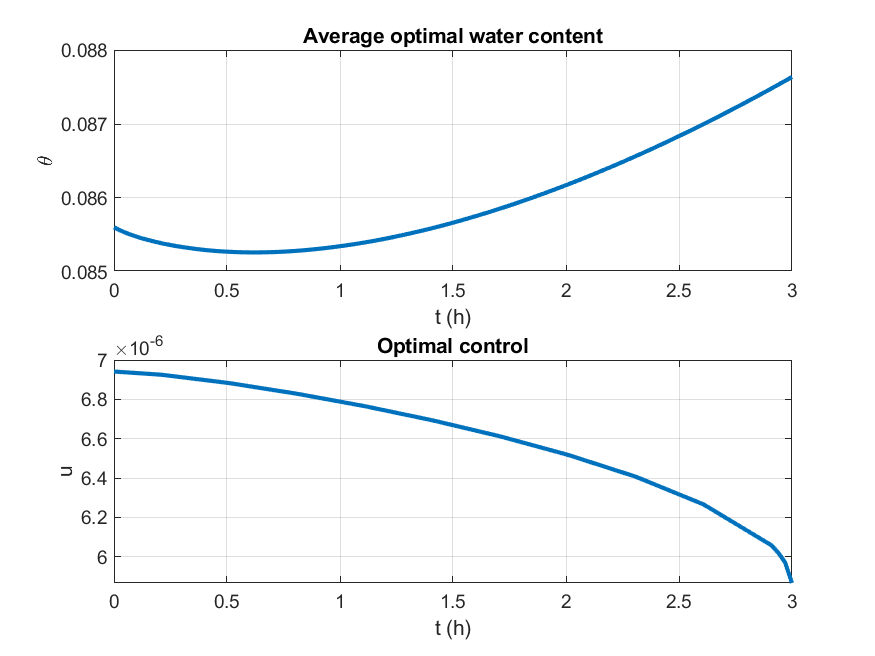}\label{fig:optControl2}}
 \caption{Numerical simulations relative to Example \ref{ex:haverkamp2}. For this simulation, initial condition is given by \eqref{eq:IC_ex2} and boundary conditions are given by \eqref{eq:bottomBC_ex2} and \eqref{eq:topBC_ex2}, respectively.
}
  \label{fig:optResults2}
 \end{figure}
 \end{exm}

 In Example \ref{ex:vg1} and Example \ref{ex:vg2} that follow, we consider the classical Van Genuchten-Mualem constitutive relations in the unsaturated zone, given by
\begin{subequations}\label{eq:VanGenutchen}
\begin{align}
\theta\left( \psi \right) &= \theta_r + \frac{\theta_S - \theta_r}{\left( 1 + \abs{\alpha \psi}^n \right)^m}, \quad m \udef 1- \frac{1}{n}, \label{eq:VanGenutchenpsi} \\
k(\psi) &= K_S  \left[ \frac{1}{1 + \abs{\alpha \psi}^n} \right]^{\frac{m}{2}} \left[ 1 - \left(1 - \frac{1}{1 + \abs{\alpha \psi}^n} \right)^m\right]^2. \label{eq:VanGenutchenk}
\end{align}
\end{subequations}
In order to verify under which conditions Van Genuchten-Mualem model satisfies the quasi-unsaturated model, we need to analyze its corresponding function $\beta(\theta(h))$.
Letting
\[
\varphi(h)\udef\frac{1}{1+\abs{\alpha h}^n},
\]
from \eqref{eq:VanGenutchen} and exploiting the fact that $h<0$, it follows that
\begin{equation}\label{eq:vgModel}
\beta(\theta(h))=\frac{K_S\left[1-(1-\varphi(h))^m\right]^2}{mn\alpha^n(\theta_S-\theta_r)\abs{h}^{n-1}\varphi(h)^{\frac{m}{2}+1}}.
\end{equation}
It is an easy computation that $\beta$ is Lipschitz, monotonically increasing and bounded from below. In order to satisfy assumption ($\mathbf{H_\beta}$), there needs
\[
\lim_{h\nearrow0^-}\beta(\theta(h))=+\infty.
\]
From~\eqref{eq:vgModel}, this condition is satisfied if and only if $n>1$.

This is the case for the simulations reported below. More specifically, we are going to consider a Berino loamy fine sand and a Glendale clay loam, with parameters drawn from \cite{Hills_et_al_1989,Berardi_Difonzo_Vurro_Lopez_ADWR_2018}.

 \begin{exm}\label{ex:vg1}
The Berino loamy fine sand is defined by the following hydraulic parameters:
\[
\theta_r = 0.0286, \theta_S = 0.3658, \alpha = 0.0280, n = 2.2390, K_S = 22.5416\,\ \mathrm{cm/h}.
\]
Here, as in Example \ref{ex:haverkamp2}, we consider a time-varying bottom condition
\begin{equation}\label{eq:bottomBC_ex3}
\theta(Z,t)=\theta_{\textrm{bottom}}(t)=\left(1-\frac{t}{T}\right)\theta_{b1}+\frac{t}{T}\theta_{b2},\,\,t\in[0,T],
\end{equation}
where
\[
\theta_{b1}\udef0.3\theta_r+0.7\theta_S,\,\theta_{b2}\udef0.1\theta_r+0.9\theta_S;
\]
boundary condition at the top of the domain is, again as in Example \ref{ex:haverkamp2},
\begin{equation}\label{eq:topBC_ex3}
\theta(0,t)=\theta_\textrm{top}(t)=u(t)+\theta_r.
\end{equation}
However, now initial condition is a quadratic polynomial function of depth and is set as
\begin{equation}\label{eq:IC_ex3}
\theta(z,0)\udef(\theta_{\rm bottom}(0)-\theta_{\rm top}(0))\left(\frac{z}{Z}\right)^2+\theta_{\rm top}(0),\quad z\in[0,Z].
\end{equation}
Results relative to this soil are depicted in Figure \ref{fig:optResults3_loamyfinesand} and are obtained using \texttt{MATLAB active-set}, where $Z=50$ cm and $T=12$ hours. \\
\end{exm}

\begin{exm}\label{ex:vg2}
Simulations on the Glendale clay loam are obtained using the following parameters:
\[
\theta_r = 0.1060, \theta_S = 0.4686, \alpha = 0.0104, n = 1.3954, K_S = 0.5458\,\, \mathrm{cm/h}.
\]
For this experiment, we fix $Z=30$ cm and $T=36$ hours; bottom boundary condition is given by
\begin{equation}\label{eq:bottomBC_ex4}
\theta(Z,t)=\theta_{\textrm{bottom}}(t)=\left(1-\frac{t}{T}\right)\theta_{b1}+\frac{t}{T}\theta_{b2},\,\,t\in[0,T],
\end{equation}
where
\[
\theta_{b1}\udef0.5\theta_r+0.5\theta_S,\,\theta_{b2}\udef0.7\theta_r+0.3\theta_S,
\]
and top boundary condition is, as in previous examples,
\begin{equation}\label{eq:topBC_ex4}
\theta(0,t)=\theta_\textrm{top}(t)=u(t)+\theta_r.
\end{equation}
Initial condition is again set as
\begin{equation}\label{eq:IC_ex4}
\theta(z,0)\udef(\theta_{\rm bottom}(0)-\theta_{\rm top}(0))\left(\frac{z}{Z}\right)^2+\theta_{\rm top}(0),\quad z\in[0,Z].
\end{equation}
Results are depicted in Figure \ref{fig:optResults3_clayloam}. Here, we employed to \texttt{MATLAB sqp} for solving the optimization problem by PGD.

\begin{figure}
 \centering
 \subfloat[][Optimal water content in the spatio-temporal domain.]
 {\includegraphics[width=.48\textwidth]{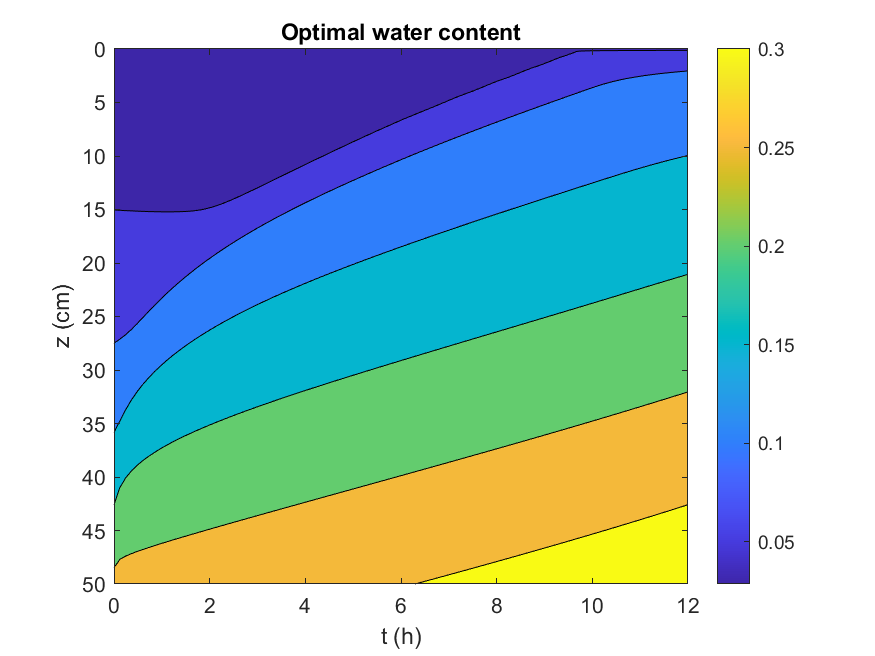}\label{fig:optTheta3}} \quad
 \subfloat[][Average optimal water content over time.]
 {\includegraphics[width=.48\textwidth]{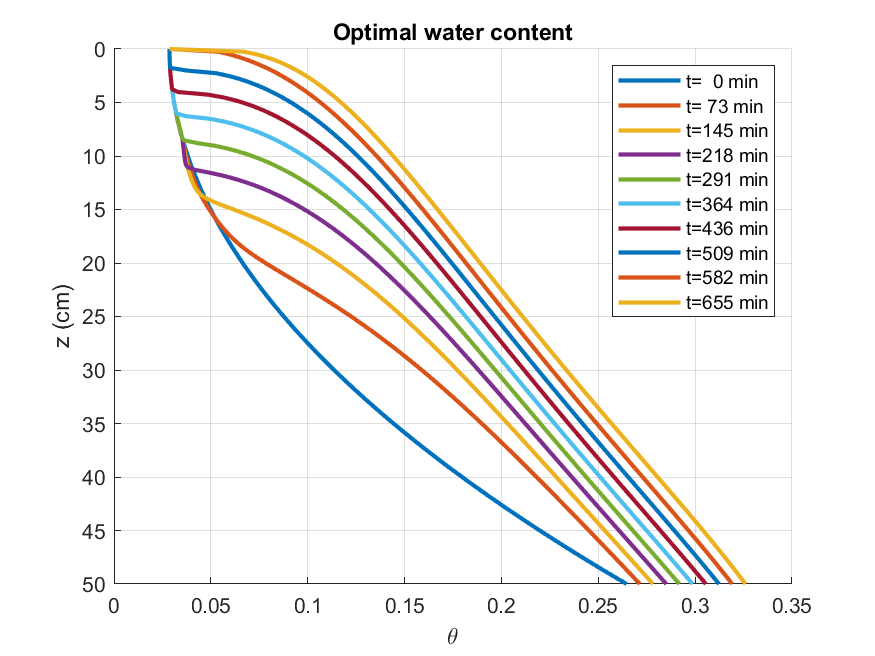}\label{fig:optMeanTheta3}} \\
 \subfloat[][Optimal water content profiles and optimal control.]
 {\includegraphics[width=.48\textwidth]{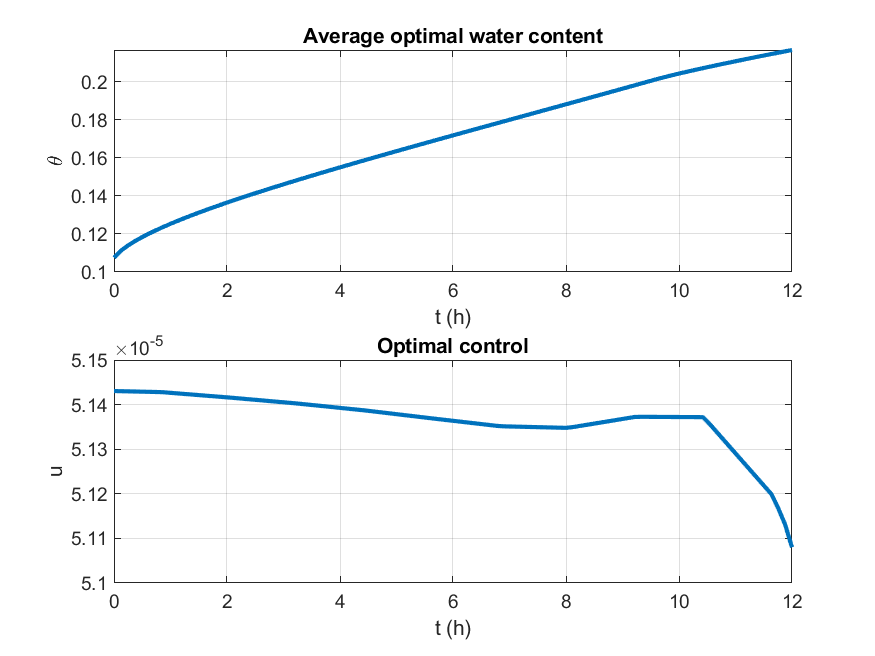}\label{fig:optControl3}}
 \caption{Numerical simulations relative to Berino loamy fine sand in Example~\ref{ex:vg1}, where we set bottom and top boundary conditions as in \eqref{eq:bottomBC_ex3} and \eqref{eq:topBC_ex3}, respectively, whilst cubic initial condition is as in \eqref{eq:IC_ex3}. 
 }
  \label{fig:optResults3_loamyfinesand}
 \end{figure}

 \begin{figure}
 \centering
 \subfloat[][Optimal water content in the spatio-temporal domain.]
 {\includegraphics[width=.48\textwidth]{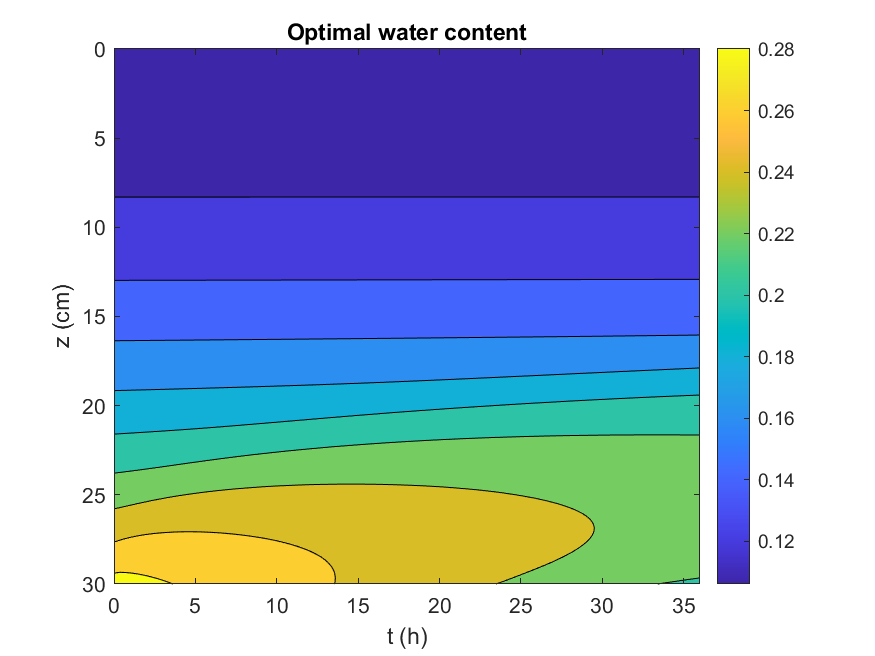}\label{fig:optTheta3_clayloam}} \quad
 \subfloat[][Average optimal water content over time.]
 {\includegraphics[width=.48\textwidth]{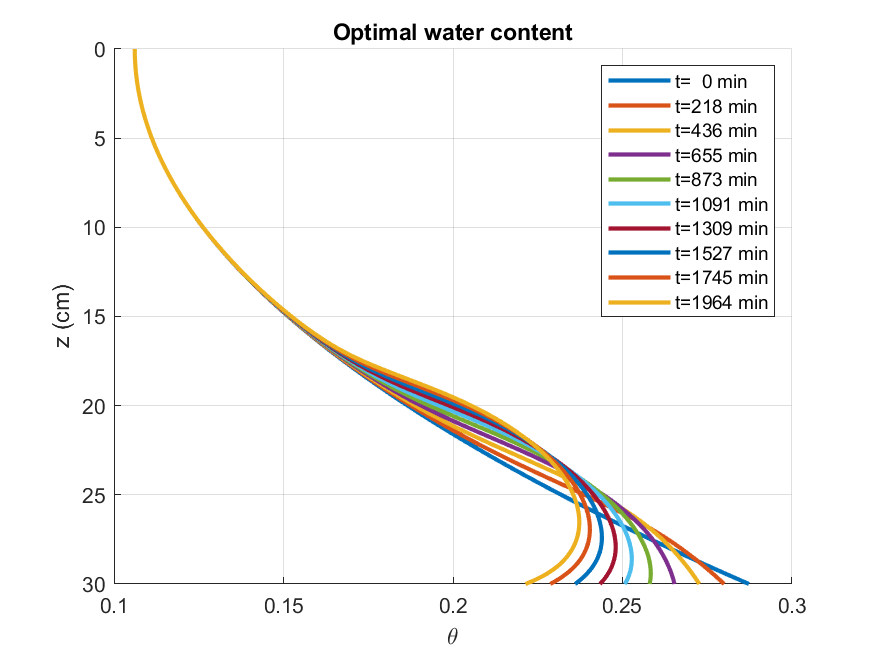}\label{fig:optMeanTheta3_clayoloam}} \\
 \subfloat[][Optimal water content profiles and optimal control.]
 {\includegraphics[width=.48\textwidth]{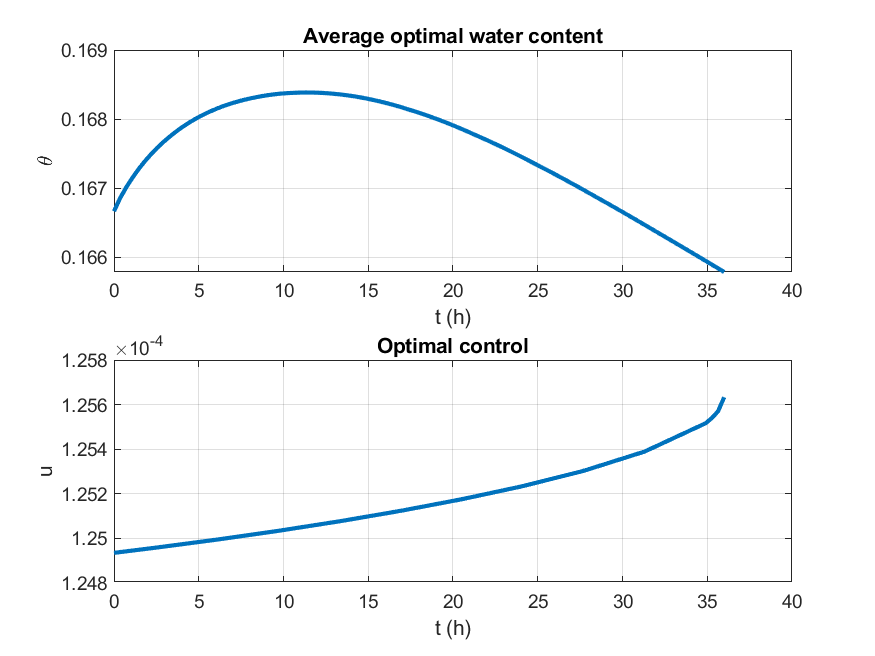}\label{fig:optControl3_clayloam}}
 \caption{Numerical simulations relative to Glendale clay loam in Example \ref{ex:vg2}. Bottom condition is given by \eqref{eq:bottomBC_ex4}, top condition by \eqref{eq:topBC_ex4}; initial condition is a quadratic polynomial as defined in \eqref{eq:IC_ex4}. 
 }
  \label{fig:optResults3_clayloam}
 \end{figure}
 \end{exm}

\section{Conclusions}

In this paper we introduce an optimal control approach aimed at optimizing the water content provided by irrigation, applying Richards' equation for unsaturated flow. We make use of quasi-unsaturated model introduced in \cite{Marinoschi_book}, extending the well-posedness results for nonlinear sink terms and deriving suitable optimality conditions for an irrigation performance index of tracking type. \\
We set the model within a {\tt MATLAB} solver by implementing a properly adapted Projected Gradient Descent method, and provide significant numerical results over a meaningful variety of soils; a deeper analytical treatise of the control system is beyond the scopes of this paper, and it is currently under investigations by the authors. \\
This paper could pave the way to an extensive use of control techniques for optimizing irrigation in real life applications, and this framework could easily be incorporated in existing irrigation software based on Richards' equation solvers. Moreover, it is worth investigating qualitative features of the more general \emph{saturated-unsaturated} model, for which there is an increasing need of both numerical and analytical results and approaches. In this context, tools from set-valued analysis and discrete control techniques could carry improvements in understanding such problems.

\section*{Acknowledgments}

MB acknowledges the partial support of RIUBSAL project funded by Regione Puglia under the call “P.S.R. Puglia 2014/2020 - Misura 16 – Cooperazione - Sottomisura 16.2 “Sostegno a progetti pilota e allo sviluppo di nuovi prodotti, pratiche, processi e tecnologie”: in particular he  thanks Mr. Giuseppe Leone and Mrs. Gina Dell'Olio for supporting the project activities; FVD has been supported by \textit{REFIN} Project, grant number 812E4967, funded by Regione Puglia: both authors acknowledge the partial support of GNCS-INdAM. RG acknowledges the support of the Natural Sciences and Engineering Research Council of Canada (NSERC), funding reference number RGPIN-2021-02632.

\bibliographystyle{elsarticle-num}
\bibliography{COMG}

\end{document}